\pdfoutput=1
\documentclass[11pt]{article}

\usepackage[utf8]{inputenc}
\usepackage[margin=1in]{geometry}
\usepackage{graphicx}
\usepackage{xcolor}
\usepackage{natbib}
\usepackage{centernot}
\usepackage{csquotes}

\usepackage{amsmath,amsfonts,bm}

\def\eqref#1{\ref{#1}}

\def\1{\bm{1}}

\DeclareMathAlphabet{\mathsfit}{\encodingdefault}{\sfdefault}{m}{sl}
\SetMathAlphabet{\mathsfit}{bold}{\encodingdefault}{\sfdefault}{bx}{n}

\newcommand{\E}{\mathbb{E}}

\newcommand{\R}{\mathbb{R}}

\newcommand{\KL}{\mathrm{KL}}

\newcommand{\Rc}{\mathcal{R}}

\newcommand{\Oc}{\mathcal{O}}

\newcommand{\Ac}{\mathcal{A}}

\newcommand{\Hc}{\mathcal{H}}
\newcommand{\Qc}{\mathcal{Q}}

\newcommand{\Pc}{\mathcal{P}}
\newcommand{\Fc}{\mathcal{F}}

\newcommand{\ito}{It\^{o}}

\renewcommand{\exp}[1]{\mathrm{e}^{#1}}

\newcommand{\dd}{\mathrm{d}}

\let\P\relax
\DeclareMathOperator{\P}{\mathbb{P}}
\DeclareMathOperator{\Q}{\mathbb{Q}}
\DeclareMathOperator{\law}{law}

\newcommand{\aux}[1]{{{#1}^{\text{aux}}}}
\newcommand{\sh}[1]{{{#1}^{\text{SHIFT}}}}

\newcommand{\pot}{U}

\newcommand{\inner}[2]{\left\langle#1,#2\right\rangle}

\newcommand{\friction}{\gamma}
\newcommand{\lip}{L}

\newcommand{\ie}{\textit{i.e.}}
\newcommand{\cf}{\textit{cf.}}

\renewcommand{\epsilon}{\varepsilon}

\newcommand{\scf}{\mathrm{sc}}
\newcommand{\rcf}{\mathrm{rc}}

\newcommand{\tc}{\hat t}
\newcommand{\ud}{\mathrm{UD}}

\usepackage{amsmath, amsthm, amssymb}
\usepackage[colorlinks=true,linkcolor=blue!60!black,citecolor=blue!60!black,urlcolor=blue!60!black]{hyperref}
\usepackage{url}
\usepackage{glossaries-extra}
\usepackage{cleveref}
\usepackage{mathtools}
\newabbreviation{lsi}{LSI}{log-Sobolev inequality}
\newabbreviation{pi}{PI}{Poincar\'e inequality}
\newabbreviation{sde}{SDE}{stochastic differential equation}
\newabbreviation{hmc}{HMC}{Hamiltonian Monte Carlo}
\newabbreviation{pde}{PDE}{partial differential equation}
\newabbreviation{ode}{ODE}{ordinary differential equation}
\newabbreviation{fpe}{FPE}{Fokker--Planck equation}
\newabbreviation{mc}{MC}{Markov chain}
\newabbreviation{mcmc}{MCMC}{Markov chain Monte Carlo}
\newabbreviation{ula}{ULA}{unadjusted Langevin algorithm}
\newabbreviation{daz}{DAZ}{diffusion at absolute zero}
\newabbreviation{smc}{SMC}{sequential Monte-Carlo}
\newabbreviation{ou}{OU}{Ornstein--Uhlenbeck}
\newabbreviation{em}{EM}{Euler-Maruyama}
\newabbreviation{angaul}{ANGAUL}{annealed gradient adjusted underdamped Langevin dynamics}
\newabbreviation{mh}{MH}{Metroplois-Hastings}
\newabbreviation{gaul}{GAUL}{gradient adjusted underdamped Langevin dynamics}
\newabbreviation{name}{ANGAUL}{annealed gradient adjusted underdamped Langevin dynamics}
\newabbreviation{pogmdm}{POGMDM}{Product of Gaussian Mixture Diffusion Model}
\newabbreviation{gmm}{GMM}{Gaussian mixture model}
\newabbreviation{kl}{KL}{Kullback--Leibler}

\newtheorem{definition}{Definition}[section] 
\newtheorem{lemma}{Lemma}[section] 
\newtheorem{theorem}{Theorem}[section] 
\newtheorem{remark}{Remark}[section]
\newtheorem{assumption}{Assumption}[section]
\crefname{assumption}{assumption}{assumptions}
\Crefname{assumption}{Assumption}{Assumptions}

\crefname{example}{example}{examples}
\Crefname{example}{Example}{Examples}
\crefname{definition}{definition}{definitions}
\Crefname{definition}{Definition}{Definitions}
\crefname{coroll}{corollary}{corollaries}
\Crefname{coroll}{Corollary}{Corollaries}
\crefname{equation}{}{}
\Crefname{equation}{Equation}{Equations}

\title{Reverse Transportation and Harnack Inequalities for Underdamped Langevin in the Non-convex Setting}

\usepackage{authblk}
\author[1,2]{Andreas Habring\thanks{Equal contribution.}}
\author[2]{Laurenz Nagler\protect\footnotemark[1]}
\affil[1]{\'Ecole polytechnique, Paris}
\affil[2]{Graz University of Technology}
\affil[ ]{\texttt{\{andreas.habring, lnagler\}@tugraz.at}}
\date{\today}

\begin{document}

\maketitle

\begin{abstract}
We establish a dimension-free long-time reverse transportation, respectively Harnack inequality for the underdamped Langevin dynamics with non-convex potential. Specifically, we prove that for any two initial points $z,z'$ and $P_t$ the underdamped Langevin semigroup, it holds
\begin{equation}
    \KL(\delta_z P_t|\delta_{z'}P_t)
    \leq C\exp{-\chi t}|z-z'|(\Rc\vee|z-z'|)
\end{equation}
for all $t>\delta$ where $\delta>0$ is arbitrarily small and $C=C(\delta),\Rc,\chi>0$ are dimension free constants. In particular, our proof extends recent results which are restricted to either the convex underdamped or the non-convex but overdamped setting.
\end{abstract}

\section{Introduction}
In this article we are interested in the underdamped Langevin \gls{sde}, that is, the stochastic process $Z_t = (X_t,V_t)$ satisfying
\begin{equation}\label{eq:uld}
    \begin{cases}
        \dd X_t =& V_t \dd t\\
        \dd V_t =& \{-\friction V_t - \nabla\pot(X_t)\}\dd t +  \sqrt{2\friction} \dd W_t\tag{ULD}
    \end{cases}
\end{equation}
with Brownian motion $(W_t)_t$, a friction parameter $\friction>0$ and a potential $\pot:\R^d\rightarrow \R$. Specifically, for possibly non-convex $\pot$ we derive a \emph{reverse transportation inequality} for~\cref{eq:uld}. That is, denoting the semigroup induced by~\cref{eq:uld} as $P_t$, we prove the existence of dimension-free constants $C,\chi, R>0$ such that 
\begin{equation}\label{eq:reverse_transpt_inequ_general}
    \begin{aligned}
        \KL(\delta_z P_t|\delta_{z'}P_t) \leq C\exp{-\chi t}|z-z'|(\Rc\vee|z-z'|),
    \end{aligned}
\end{equation}
where $\KL$ denotes the Kullback--Leibler divergence. We may represent $\delta_zP_t$ and $\delta_{z'}P_t$ as the laws of $Z_t$ respectively $Z'_t$ which solve \eqref{eq:uld} with initial condition $Z_0 = z$ , $Z'_0 = z'$ and $z,z' \in \R^{2d}$, that is, 
\begin{equation}\label{eq:intro_coupling}
    \begin{aligned}
        \begin{cases}
            \dd X_t &= V_t \dd t\\
            \dd V_t &= \{-\friction V_t - \nabla\pot(X_t)\}\dd t +  \sqrt{2\friction} \dd W_t,\\
            (X_0,V_0) &= z
        \end{cases}\\
        \begin{cases}
            \dd X'_t &= V'_t \dd t\\
            \dd V'_t &= \{-\friction V'_t - \nabla\pot(X'_t)\}\dd t +  \sqrt{2\friction} \dd W_t\\
            (X'_0,V'_0) &= z'.
        \end{cases}
    \end{aligned}
\end{equation}
Moreover, \cref{eq:reverse_transpt_inequ_general} is equivalent to a log-Harnack inequality, that is, for any sufficiently well-behaved $f:\R^{2d}\rightarrow (0, \infty)$
\begin{equation}\label{eq:harnack_general}
    P_t(\log f)(z) \leq \log P_t f (z') +  C\exp{-\chi t}|z-z'|(R\vee|z-z'|).
\end{equation}
We refer to~\cref{sec:harnack} for a short proof via Donsker-Varadhan variational formula.

A main obstruction in proving~\cref{eq:reverse_transpt_inequ_general} is the fact that the initial distributions $\delta_z, \delta_{z'}$ are mutually singular whenever $z\neq z'$ prohibiting many arguments based on a direct application of Girsanov's theorem to compare the two \glspl{sde} in~\cref{eq:intro_coupling}.

As a remedy, one may consider coupling techniques using a third auxiliary process, a technique referred to as \emph{shifted Girsanov}, \cf~\citep{lu2026longtimereversetransportationinequalities,altschuler2025shifted_comp_1,altschuler2026shiftedcompositionivballistic,arnaudon2006harnackinequalityheatkernel}. More specifically, we consider an auxiliary bridge process \((Z^\aux{}_t)_{t}\) interpolating between $z'$ and $Z_T$, \ie, $Z^\aux{}_0=z'$ and \(Z^\aux{}_T = Z_T\) at a prescribed time \(T > 0\). This allows us to write, \(\KL(Z_T\mid Z_T') = \KL(Z^\aux{}_T\mid Z_T')\). By designing the process $(Z^\aux{}_t)_t$ appropriately, we can ensure that the path measures of $(Z^\aux{}_t)_t$ and $(Z_t')_t$ are, indeed, mutually absolutely continuous hence allowing to estimate $\KL(Z^\aux{}_T\mid Z_T')$ using Girsanov.

\section{Related Work and 
Contributions}\label{sec:related_work}
\paragraph{Underdamped Langevin} 
Due to its accelerating properties there is a vast literature on underdamped Langevin dynamics. The main technical difficulty of underdamped Langevin is the degeneracy of the Brownian motion, which led to Villani's development of hypocoercivity~\citep{villani2006hypocoercivity}.
In the machine learning community the research on underdamped Langevin was largely incited by~\citep{cheng2018underdamped}. For works considering the convex setting and synchronous coupling, see~\citep{Dalalyan2018OnSF,cheng2018underdamped,leimkuhler2024contraction}. In the non-convex setting, reflection coupling with carefully designed distance functions has been fruitfully used to derive contraction results in Wasserstein distances~\citep{schuh2024globalcontractivity,eberle2017coupling_kinetic,eberle2026non_rev_lift}. Furthermore, there is a vast literature on different discretizations of underdamped Langevin dynamics. Most notably, its structure is amenable to splitting methods where individual components of the drift are integrated alternatingly~\citep{leimkuhler2013rational,leimkuhler2024contraction}. A discretization scheme that has been particularly successful is the randomized midpoint scheme~\citep{altschuler2026shiftedcompositionivballistic,shen2019randomized}. Of all the cited works, the most closely related is undoubtedly~\citep{altschuler2026shiftedcompositionivballistic}. 
There the authors employ the shifted Girsanov argument that is also used in the present article and leads to improved scaling in the condition number $\kappa$. 

\paragraph{Shifted Girsanov and reverse-Transport- and Harnack Inequalities} 
Dimension-free Harnack inequalities were first shown for heat semigroups on manifolds via semigroup methods in~\citep{wang1997logarithmic} and afterward via coupling techniques in~\citep{arnaudon2006harnackinequalityheatkernel}. By now there is an extensive body of literature regarding Harnack inequalities and their applications and we refer to~\citep{wang2010couplingapplications,wang2013harnack,wang2010harnack,wang2011harnack} for some overviews.
In fact,~\citep{arnaudon2006harnackinequalityheatkernel} is the work pioneering the idea of using a third auxiliary process together with Girsanov's theorem to derive Harnack type inequalities.
More recently, the works of 
\citep{altschuler2025shifted_comp_1} and  \citep{altschuler2026shiftedcompositionivballistic} used a shifted Girsanov arguments to establish reverse transport inequalities for the overdamped and underdamped Langevin dynamics. In \citep{lu2026longtimereversetransportationinequalities} the techniques from~\citep{altschuler2026shiftedcompositionivballistic} are applied to the overdamped Langevin diffusion in the non-logconcave setting. Finally, \citep{guillin2012degeneratefokkerplanckequations} derived reverse transportation inequalities for the underdamped Langevin diffusion via a chain of implications (\cf~\citep[Section 1.3]{wang2013harnack}). However, in contrast to our results, their bounds do not converge to $0$ as $T \rightarrow \infty$.

\paragraph{Contributions}
As mentioned above, our work can be seen as an extension of some of the results in~\citep{altschuler2026shiftedcompositionivballistic,lu2026longtimereversetransportationinequalities}. In particular,~\citep{altschuler2026shiftedcompositionivballistic} provide a corresponding convergence result only for convex potentials\footnote{Note that the authors obtain the Harnack inequality more generally. However, the convergence to zero is only shown in the convex case.}
and~\citep{lu2026longtimereversetransportationinequalities} for non-convex potentials but for the overdamped Langevin dynamics. Conversely the results in the present article cover the non-convex and underdamped setting. It should be mentioned, however, that~\citep{altschuler2026shiftedcompositionivballistic,lu2026longtimereversetransportationinequalities} provide additional results in other contexts. In particular, in~\citep{altschuler2026shiftedcompositionivballistic} the authors also consider the discrete setting and in~\citep{lu2026longtimereversetransportationinequalities} results are also obtained for the, more general, R\'enyi divergence.


\section{Notation and Assumptions}\label{sec:Notation_Assumption}

\paragraph{Notation}
Throughout $x$ will be used for the position variable, $v$ for the velocity, and $z=(x,v)$ for the position-velocity pair and we note that several different variations of these symbols will appear ($x$, $X$, $X'$, $X^\aux{}$ etc.).
The set of probability measures on \(\R^d\) is denoted as \(\Pc(\R^d)\). We write \(\mu \ll \nu\) when \(\mu\) is absolutely continuous w.r.t. \(\nu\), and with \(\frac{\dd \mu}{\dd \nu}\) for the Radon-Nikodym derivative of \(\mu\) w.r.t. \(\nu\). We define the Kullback-Leibler divergence as 
\begin{equation}
    \begin{aligned}
        \KL: \Pc(\R^d)\times \Pc(\R^d)&\rightarrow [0,\infty]\\
        (\mu,\nu)&\mapsto\KL(\mu|\nu) = \begin{cases}
            \int \log(\frac{\dd \mu}{\dd \nu})\frac{\dd \mu}{\dd \nu} \, \dd \nu\quad &\text{if $\mu\ll \nu$}\\
            \infty\quad &\text{else.}
        \end{cases}
    \end{aligned}
\end{equation}
We will denote for simplicity for two random variables $X,Y$
\begin{equation}
    \KL(X|Y) \coloneqq \KL(\law(X)|\law(Y)).
\end{equation}
The symbol $\tilde{\Oc}$ is used to denote complexity results omitting logarithmic factors. 

Our main assumptions about the potential $\pot$ are $L$-smoothness and strong convexity outside a ball:
\begin{assumption}\label{ass:pot}
    The potential $\pot\in C^2(\R^d)$ and there exist $\lip,m,R>0$ such that we have 
    \begin{align}
        \lip I\succeq\nabla^2 \pot(x)\succeq& -\lip I,\quad x\in \R^d\label{ass:lip}\\
        \nabla^2 \pot(x)\succeq& mI,\quad |x|\geq R.\label{ass:convex}
    \end{align}
\end{assumption}

\section{Proof of the main result}\label{sec:cont_time_analysis}

\subsection{Set-up}
To prove the above stated reverse transportation inequality we rely on the \emph{shifted composition} principle introduced in~\citep{arnaudon2006harnackinequalityheatkernel} and recently \emph{rediscovered} in~\citep{altschuler2025shifted_comp_1,altschuler2026shiftedcompositionivballistic}. We explain and motivate the approach in the following. 

Our goal is to bound the $\KL$ divergence between the law of two solutions of~\cref{eq:uld} initialized at $z_0,z'_0 \in \R^{2d}$, respectively
\begin{equation}
    \begin{aligned}
        &\begin{cases}
            \dd X_t &= V_t \dd t\\
            \dd V_t &= \{-\friction V_t - \nabla \pot ( X_t) \}\,\dd t + \sqrt{2 \friction}\dd W_t\\
            Z_0 = z_0
        \end{cases}\\
        &\begin{cases}
            \dd X_t' &= V'_t \dd t\\
            \dd V_t' &= \{-\friction V_t' - \nabla \pot (X_t') \}\,\dd t + \sqrt{2 \friction}\dd W_t\\
            Z'_0 = z_0'.
        \end{cases}
    \end{aligned}
\end{equation}
Recall, however, that a na\"ive application of Girsanov's theorem to these processes as-is will not lead to meaningful bounds as the different initializations render the respective path measures mutually singular. As a remedy, we introduce a deliberately designed interpolating auxiliary process $Z^\aux{}$ which satisfies on one hand $Z^\aux{}_0 = z'_0$ so that we may apply Girsanov's theorem to compare $Z^\aux{}$ and $Z'$. On the other hand, we design $Z^\aux{}$ such that at a finite time $T>0$ it holds $\law(Z^\aux{}_T)=\law(Z_T)$ which implies $\KL(Z_T|Z'_T) = \KL(Z^\aux{}_T|Z'_T)$. In other words, the process $Z^\aux{}$ constitutes a bridge from $Z'_0$ to $Z_T$. Let us now make this construction more concrete:

For \(z_0, z_0' \in \R^{2d}\) we define the coupled \gls{sde},
\begin{equation}\label{eq:coupled_uld}
    \begin{aligned}
    &\begin{cases}
        \dd X_t &= V_t \dd t\\
        \dd V_t &= \{-\friction V_t - \nabla \pot (X_t) \}\,\dd t + \sqrt{2 \friction}\dd W_t\\
        \dd X_t^{\aux{}} &= V_t^{\aux{}}\dd t\\
        \dd V_t^{\aux{}} &=
        \begin{aligned}[t]
             &\{-\friction V_t^{\aux{}}
            - \nabla \pot (X^\aux{}_t)\\
            & + \eta_t^x(X_t - X_t^{\aux{}})
            + \eta_t^v(V_t - V_t^{\aux{}})\}\dd t
            + \sqrt{2 \friction}\,\dd \bar W_t.
        \end{aligned}\\
        Z_0 &= z_0\\
        Z_0^\aux{} &= z_0'\\
    \end{cases}
    \end{aligned}\tag{AUX}
\end{equation}
where $(W_t)_t$ and $(\bar{W}_t)_t$ are two Brownian motions coupled in a deliberate way, which will be explained below. Crucially, the shift parameters $\eta^v_t$, $\eta^x_t$ will be chosen to satisfy $\eta^v_t,\eta^x_t\rightarrow \infty$ as $t \rightarrow T$ ensuring that $Z^\aux{}_T = Z_T$ at a finite time $T$.
In addition, we have a second solution of~\cref{eq:uld} initialized at $z_0'$
\begin{equation}\label{eq:girs_uld}
    \begin{aligned}
    &\begin{cases}
        \dd X_t' &= V_t'\dd t \\
        \dd V_t' &=
        \begin{aligned}[t]
             &\{-\friction V_t'
            - \nabla \pot(X_t')\}\dd t
            + \sqrt{2 \friction}\,\dd W_t'.
        \end{aligned}\\
        Z'_0 &= z_0'
    \end{cases}
    \end{aligned}
\end{equation}
where we note that, technically,~\cref{eq:girs_uld} is not defined on the same probability space as~\cref{eq:coupled_uld} (see~\cref{sec:rigorous_girsanov} for details).
Applying Girsanov's theorem to compare $Z^{\aux{}}$ with $Z'$ then yields
\begin{equation}\label{eq:girsanov_int_full_time}
    \begin{aligned}
        \KL(Z_T|Z_T')
        = \KL(Z^{\aux{}}_T|Z'_T)
        &\leq \KL((Z^{\aux{}}_t)_{t\in [0,T]}|(Z'_t)_{t\in [0,T]})\\
        &= \frac{1}{4\friction}\int_0^T \E[|\eta_t^x(X_t - X_t^{\aux{}})
                + \eta_t^v(V_t - V_t^{\aux{}})|^2]\dd t
    \end{aligned}
\end{equation}
where we used the hitting condition $\law(Z^\aux{}_T)=\law(Z_T)$ and the data processing inequality. 
To now obtain the desired bounds, \cref{eq:reverse_transpt_inequ_general} respectively \cref{eq:harnack_general}, it remains to estimate the above integral. 

\paragraph{Proof outline}
The proof consists of the following steps:
\begin{enumerate}
    \item Inspired by~\citep{lu2026longtimereversetransportationinequalities} we define a time $\tc\in [0,T]$ at which a \emph{regime transition} happens. 
    For $t\geq \tc$ we will choose $\eta^x_t,\eta^v_t\gg 0$. In this regime the shifts are strong, and a synchronous coupling, that is, $\dd W_t = \dd \bar W_t$ in~\cref{eq:coupled_uld} is sufficient to show contraction and prove the hitting condition $Z^\aux{}_T=Z_T$.
    \item Conversely, for $t < \tc$ we set $\eta^x_t = \eta^v_t = 0$. In this regime, $Z^\aux{}_t$ and $Z_t$ are two conventional \eqref{eq:uld} processes. Thus, we employ a reflection coupling for $W_t$ and $\bar W_t$ (\cf~\citep{eberle2017coupling_kinetic,schuh2024globalcontractivity}) which allows us to derive contraction in a deliberately designed distance without any shifts.
    \item The result follows by concatenating the bounds in the two regimes.
\end{enumerate}
The remainder of the article contains the proof details.
\subsection{Coupling construction}\label{sec:coupling}
It is well-known in the literature that purely synchronous coupling is unsuitable to derive contraction results for~\cref{eq:uld} in the non-convex setting (\cf~\citep{eberle2017coupling_kinetic}). Thus, we rely on a combination of synchronous and reflection coupling. More precisely, building on~\citep{eberle2017coupling_kinetic,schuh2024globalcontractivity} we first introduce new coordinates. Specifically, let us denote
\begin{equation}\label{eq:twisted_coordinate1}
    P_t = X_t + \friction^{-1}V_t
\end{equation}
and analogously for the auxiliary process. Then~\cref{eq:uld} reads as
\begin{equation}\label{eq:uld_twisted}
    \begin{cases}
        \dd X_t =& -\friction X_t\dd t + \friction P_t \dd t\\
        \dd P_t =& -\friction^{-1}\nabla \pot(X_t)\dd t + \sqrt{2\friction^{-1}}\dd W_t.
    \end{cases}
\end{equation}
Now we define the processes $(W_t)_t$ and $(\bar W_t)_t$ in~\cref{eq:coupled_uld} using a combination of synchronous and reflection coupling as
\begin{align}\label{eq:reflection_brownian_motions}
    W_t &= \int_0^t \rcf(Z_s, Z_s^\aux{}) \dd W_s^{\rcf} + \int_0^t \scf(Z_s, Z_s^\aux{})\dd W_s^{\scf},\\
    \bar W_t &= 
    \begin{cases}
        \int_0^t \rcf(Z_s, Z_s^\aux{})R_s \dd W_s^{\rcf} + \int_0^t \scf(Z_s, Z_s^\aux{})\dd W_s^{\scf} \quad \text{for \(t < \tc\)} \\
        W_t - W_{\tc} + \bar  W_{\tc} \quad \text{for \(t \geq \tc\)}
    \end{cases}
\end{align}
with $(W^{\rcf}_t)_t$, $(W^{\scf}_t)_t$ two independent Brownian motions and the reflector \(R_s = (I - 2e_s e_s^\top)\) with \(e_t = Q_t/|Q_t|\) and $Q_t = P_t-P^\aux{}_t$.
The two smooth cut-off functions \(\scf, \rcf: \R^{2d} \times \R^{2d} \rightarrow [0,1]\) satisfy \(\scf^2 + \rcf^2 = 1\). One may then check that \(W_t\) and \(\bar W_t\) are, indeed,  Brownian motions by Levy's characterization~\cite[Theorem 3.16]{karatzas2014brownian}.
While for $t>\tc$ we have standard synchronous coupling, the construction for $t<\tc$ is exactly the one from~\citep{schuh2024globalcontractivity} and we refer to this work for the details and limit ourselves to a high-level explanation here.

Essentially, the coupling should be synchronous whenever $Q_t=0$ or $Z_t-Z^\aux{}_t$ is large. The reason therefore is that in the first case, $P_t=P^\aux{}_t$ and thus the first equation in~\cref{eq:uld_twisted} yields a contraction also in the $X$-component. In the second case, on the other hand, the convexity outside a ball of $U$ provides contraction with synchronous coupling. Thus, the reflection coupling is active only in the intermediate regime, where $Z_t-Z^\aux{}_t$ is small but not zero.


\subsection{Contraction for $t>\tc$}
We begin by investigating the regime $t>\tc$. In this setting, we define a slightly different twisted coordinate system than~\cref{eq:twisted_coordinate1}. Specifically, we set
\begin{equation}
    \begin{aligned}
        \begin{bmatrix}
            X_t\\
            \bar P_t
        \end{bmatrix} \coloneqq 
        \begin{bmatrix}
            X_t\\
            X_t + \frac{2}{\friction_t}V_t
        \end{bmatrix}
        = A_t \begin{bmatrix}
            X_t\\
            V_t
        \end{bmatrix},\quad 
        A_t \coloneqq \begin{bmatrix}
            1&0\\
            1&\frac{2}{\friction_t}
        \end{bmatrix},\quad \friction_t = \friction + \eta^v_t.
    \end{aligned}
\end{equation}
The definition of the time-dependent friction parameter $\friction_t = \friction + \eta^v_t$ follows naturally from collecting the $V^\aux{}_t$-terms in~\cref{eq:coupled_uld}. Note that $\bar P$ differs from $P$ by the factor $2$ in the bottom right entry of $A_t$.
Moreover, we introduce the notation
\begin{equation}
    D_t \coloneqq \begin{bmatrix}
            X_t-X_t^\aux{}\\
            V_t - V^\aux{}_t
        \end{bmatrix},\quad T_t \coloneqq A_t D_t
\end{equation}
for the difference processes in regular and twisted coordinates.
Next, we fix the precise definition of the shift parameters.
\begin{definition}[Shift parameters]
We define the shift parameters for $c_0>0$, as
    \begin{equation}
    \eta^v_t = \1_{t\geq\tc}\frac{c_0\omega}{\exp{\omega(T-t)}-1},\quad \eta^x_t = \frac{\gamma_t\eta^v_t}{2}
\end{equation}
where
\begin{equation}
     \omega =\begin{cases}
         -\frac{\lip}{3\friction}\quad &\text{if $\friction\geq \sqrt{32\lip}$ (high friction)}\\
         -\frac{\sqrt{\lip}}{3}\quad &\text{if $\friction< \sqrt{32\lip}$ (low friction)}
     \end{cases}
\end{equation}
\end{definition}
Note that, indeed, $\eta^v_t,\eta^x_t\rightarrow \infty$ as $t\rightarrow T$. The following result is shown in~\citep{altschuler2026shiftedcompositionivballistic}. We include the proof for completeness. 
\begin{lemma}\label{lemma:sync_contraction}
    There exist absolute constants $c',c_0>0$ such that for any $s,t\in [\tc, T]$, $s<t$ it holds
    \begin{equation}
        \E[|T_t|^2] \leq \bigg(\frac{1-\exp{-\omega(T-t)}}{1-\exp{-\omega(T-s)}}\bigg)^{c'c_0}
            \E[|T_s|^2].
    \end{equation}
    Moreover, choosing $c_0>0$ sufficiently large it holds $\E[|Z_t-Z^\aux{}_t|^2]\rightarrow 0$ as $t\to T$.
\end{lemma}
\begin{proof}
    Recall that for $t>\tc$ we have a synchronous coupling (\cf~\cref{eq:reflection_brownian_motions}) and it holds
    \begin{align} \label{eq:ito_sq_T_t1}
        \dd |T_t|^2 = 2\langle T_t,  \dd T_t \rangle
    \end{align}
    where, since $T_t = A_t D_t$,
    \begin{align} \label{eq:difference_sde_transformed}
        \dd T_t &= \dot A_t D_t \dd t  + A_t \dd D_t.
    \end{align}
    For the second term we find,
    \begin{align}
        \dd D_t &= 
        \begin{bmatrix}
            0 & I \\
            -\eta_t^x I - \Hc_t & - \friction_t I  
        \end{bmatrix}
        D_t \dd t
        = M_t D_t \dd t
    \end{align}
    where we used that by the fundamental theorem of calculus
    \begin{equation}
        \nabla_x \pot(X_t) - \nabla_x \pot(X_t^\aux{}) = \underbrace{\int_0^1\nabla^2 \pot(X_t^\aux{} + \lambda (X_t-X^\aux{}_t) \dd \lambda}_{\eqqcolon \Hc_t}\, (X_t - X_t^\aux{}).
    \end{equation}
    Inserting in \cref{eq:ito_sq_T_t1} yields
    \begin{align*}
        \dd |T_t|^2 
        =& 2 \inner{T_t}{(\dot A_t + A_t M_t)A_t^{-1} T_t} \dd t = -2\inner{T_t}{\Qc_t T_t} \dd t
    \end{align*}
    where we introduce the notation $\Qc_t = -\mathrm{sym}((\dot A_t + A_t M_t)A_t^{-1})$ with the symmetrization $\mathrm{sym}(M) \coloneqq (M+M^\top)/2$ for any quadratic matrix $M$.
    By~\cite[Lemma 3.13]{altschuler2026shiftedcompositionivballistic} for $c_0\geq 24$, it holds for the smallest eigenvalue of $\Qc$
    \begin{equation}
        \lambda_{\min}(\Qc_t)\geq \frac{\eta^v_t}{48}
    \end{equation}
    where we note that this bound is proven in~\citep[Lemma 3.13]{altschuler2026shiftedcompositionivballistic} in both the high and low friction settings.
    Inserting the lower bound and the definition of $\eta^v_t$ yields with $c'=\frac{1}{48}$
    \begin{equation}
        \dd  |T_t|^2 
        \leq -\frac{c'c_0\omega}{\exp{\omega(T-t)}-1}|T_t|^2\dd t.
    \end{equation}
    By Grönwall it follows for any $\tc\leq s<t$
    \begin{equation}
        \begin{aligned}
            \E[|T_t|^2]
            \leq& \left(\frac{1-\exp{-\omega(T-t)}}{1-\exp{-\omega(T-s)}}\right)^{c'c_0}
            \E[|T_s|^2].
        \end{aligned}
    \end{equation}
    To proof the hitting condition, note that since $\E[|X_t-X^\aux{}_t|^2]\leq \E[|T_t|^2]\to 0$ as $t\to T$ we directly obtain convergence of the $x$-component. For the $v$-component we have for $c'c_0>2$
    \begin{equation}
        \begin{aligned}
            \E[|V_t-V^\aux{}_t|^2]
            \lesssim \friction_t^2\E[|T_t|^2]
            \lesssim& \frac{1}{(\exp{\omega (T-t)}-1)^2}\left(\frac{1-\exp{-\omega(T-t)}}{1-\exp{-\omega(T-{\tc})}}\right)^{c'c_0}\E[|T_{\tc}|^2]\\
            \lesssim& \frac{1}{\exp{2\omega (T-t)}}\frac{(1-\exp{-\omega(T-t)})^{c'c_0-2}}{(1-\exp{-\omega(T-{\tc})})^{c'c_0}}
            \E[|T_{\tc}|^2]\to 0
        \end{aligned}
    \end{equation}
    as $t\to T$.
\end{proof}

As a consequence we obtain the following intermediate upper bound for the $\KL$ divergence.
\begin{lemma}\label{lemma:interm_bound}
It holds
\begin{equation}\label{eq:interm_bound}
    \begin{aligned}
            \KL(Z_T|Z'_T)
            \lesssim \bigg(\frac{|\omega|^3(\exp{-4\omega(T-\tc)}-1)}{(1-\exp{-\omega(T-\tc)})^4} + \frac{|\omega|\friction^2(\exp{-2\omega(T-\tc)}-1)}{(1-\exp{-\omega(T-\tc)})^2}\bigg)\E[|T_{\tc}|^2].
        \end{aligned}
    \end{equation}
\end{lemma}
\begin{proof}
    By~\cref{lemma:sync_contraction} $Z_t-Z^\aux{}_t\to 0$ as $t\to T$ in $L^2(\P)$. By $L^2(\P)$-continuity of $(Z_t)_t$ it follows
    \begin{equation}
        \E[|Z_T - Z^\aux{}_t|^2]^{1/2}\leq \E[|Z_T - Z_t|^2]^{1/2} + \E[|Z_t - Z^\aux{}_t|^2]^{1/2}\to 0.
    \end{equation}
    Thus, it holds $Z^\aux{}_t\to Z_T$ and $Z_t'\to Z_T'$ both in $L^2$. Indeed, $L^2$ convergence implies weak convergence of the laws and therefore, by (joint) weak lower semicontinuity of the $\KL$ we do, in fact, have
    \begin{equation}\label{eq:kl_bound_girsanov}
        \KL(Z_T|Z_T')
        \leq \liminf_{t\to T}\KL(Z^{\aux{}}_t|Z'_t)
        \leq \frac{1}{4\friction}\int_0^T \E[|\eta_t^x(X_t - X_t^{\aux{}})+ \eta_t^v(V_t - V_t^{\aux{}})|^2]\dd t
    \end{equation}
    where we used Girsanov's theorem in the last inequality.    
    Using~\cref{lemma:sync_contraction} we obtain
    \begin{equation}
        \begin{aligned}
            \KL(Z_T|Z_T')
            \leq& \frac{1}{4\friction}\int_{\tc}^T (\eta^v_t)^2\friction_t^2\E[|T_t|^2]\dd t\\
            \leq& \frac{1}{4\friction}\int_{\tc}^T (\eta^v_t)^2\friction_t^2 \left(\frac{1-\exp{-\omega(T-t)}}{1-\exp{-\omega(T-\tc)}}\right)^{c'c_0}\dd t\,\E[|T_{\tc}|^2]\\
        \end{aligned}
    \end{equation}
    The integral can be estimated as follows,
    \begin{align}
        \int_{\tc}^T &(\eta^v_t)^2\friction_t^2 \left(\frac{1-\exp{-\omega(T-t)}}{1-\exp{-\omega(T-\tc)}}\right)^{c'c_0}\dd t\\
        \lesssim& \int_{\tc}^T \frac{\omega^4}{(\exp{\omega(T-s)}-1)^4}\bigg(\frac{1-\exp{-\omega(T-s)}}{1-\exp{-\omega(T-\tc)}}\bigg)^{c'c_0}\dd s
        + \int_{\tc}^T \frac{\omega^2\friction^2}{(\exp{\omega(T-s)}-1)^2}\bigg(\frac{1-\exp{-\omega(T-s)}}{1-\exp{-\omega(T-\tc)}}\bigg)^{c'c_0}\dd s\\
        \intertext{noting that $\frac{1-\exp{-\omega(T-s)}}{1-\exp{-\omega(T-\tc)}}\in [0,1]$ and assuming $c_0>0$ large enough so that $c'c_0\geq 4$}
        \lesssim& \int_{\tc}^T \frac{\omega^4}{(\exp{\omega(T-s)}-1)^4}\bigg(\frac{1-\exp{-\omega(T-s)}}{1-\exp{-\omega(T-\tc)}}\bigg)^{4}\dd s
        + \int_{\tc}^T \frac{\omega^2\friction^2}{(\exp{\omega(T-s)}-1)^2}\bigg(\frac{1-\exp{-\omega(T-s)}}{1-\exp{-\omega(T-\tc)}}\bigg)^{2}\dd s\\
        \intertext{using that $\frac{1-1/u}{u-1} = \frac{1}{u}$}
        \lesssim& \int_{\tc}^T \frac{\omega^4}{\exp{4\omega(T-s)}}\frac{1}{(1-\exp{-\omega(T-\tc)})^4}\dd s + \int_{\tc}^T \frac{1}{\exp{2\omega(T-s)}}\frac{\omega^2\friction^2}{(1-\exp{-\omega(T-\tc)})^2}\dd s\\
        \lesssim& \frac{|\omega|^3(\exp{-4\omega(T-\tc)}-1)}{(1-\exp{-\omega(T-\tc)})^4} + \frac{|\omega|\friction^2(\exp{-2\omega(T-\tc)}-1)}{(1-\exp{-\omega(T-\tc)})^2}\label{eq:final_bound3}
    \end{align}
    yielding the desired result.
\end{proof}
\begin{remark}
    Note that, if we were to set $\tc=0$ the bound in~\cref{eq:interm_bound} would not tend to zero as $T\to\infty$ since $\omega<0$. The exponential decay therefore has to be established in $[0,\tc]$.
\end{remark}
\begin{remark}
    Note that the above arguments -- especially the application of Girsanov to bound the $\KL$ -- rely on strong existence and uniqueness of solutions of the occurring \glspl{sde} which is not trivially given due to the singular coefficients $\eta^v,\eta^x$. However, for any $t<T$, strong existence and uniqueness on $[0,t]$ follow from standard results due to Lipschitz-continuity. While $\eta^v_s, \eta^x_s$ exhibit a discontinuity at $s=\tc$, we have strong existence and uniqueness, on $[0,\tc)$ as well as $[t,\tc]$ leading to a unique solution on $[0,t]$ by \emph{concatenation}. This is sufficient for the limit argument in~\cref{eq:kl_bound_girsanov}.
\end{remark}

\subsection{Contraction for $t<\tc$}
In view of~\cref{eq:interm_bound} it remains to bound $\E[|T_{\tc}|^2]$. Given that $\eta^v_t=\eta^x_t=0$ for $t\in[0,\tc)$ this means we have to obtain a contraction result for \emph{regular}~\cref{eq:uld}. Readily available results, however, do not directly yield satisfying rates. Indeed, we can identify three groups of contraction results:
\begin{enumerate}
    \item Strongly convex: For instance~\citep{cheng2018underdamped,Dalalyan2018OnSF} provide contraction in the $L^2$-norm under the assumption of a strongly logconcave target using synchronous coupling. However, since we consider the non-logconcave setting, we cannot apply these results.
    \item Non-logconcave: In~\citep{eberle2017coupling_kinetic} exponential decay in $L^2$ is established using reflection and synchronous coupling even in the non-logconcave case. However, the convergence deteriorates exponentially in the dimension.
    \item Lastly, using similar techniques as in~\citep{eberle2017coupling_kinetic}, in~\citep{schuh2024globalcontractivity} dimension-free exponential decay was proven in the $L^1$-norm.
\end{enumerate}
Our goal is to apply the above mentioned $L^1$-convergence result from~\citep{schuh2024globalcontractivity} which we state explicitly in the following.
\begin{lemma}[Adapted from Theorem 5 and Remark 3 in \citep{schuh2024globalcontractivity}]\label{lemma:schuh_contraction}
    Let~\cref{ass:pot} hold.
    Let further \(\friction^2 > \frac{2(L+m)^2}{m}\) then it holds for some dimension free constants $\chi,  C>0$ and all $0\leq s\leq t<\tc$
    \begin{equation}
        \E_{\P}[|T_t|] \leq  C\exp{- \chi (t-s)}\E_{\P}[|T_s|].
    \end{equation}
\end{lemma}
In order to apply~\cref{lemma:schuh_contraction} we have to go from the $L^2$ to the $L^1$ norm.
By Hölder's inequality we have
\begin{equation}
    \E[|T_{\tc}|^2]\leq \E[|T_{\tc}|] \;\|T_{\tc}\|_{L^\infty(\P)}.
\end{equation}
Thus, to conclude, we are left to derive an a.s. uniform bound of $T_{\tc}$. 

\begin{lemma}\label{lemma:uniform_bound}
    Let~\cref{ass:pot} hold and assume $\friction\geq \sqrt{2\lip}$. Then, for $t\in [0,\tc]$ a.s.
    \begin{equation}
        |T_t|\lesssim (1\vee \friction^{-1}) (\Rc\vee|z_0-z_0'|).
    \end{equation}
    where
    \begin{equation}
        \Rc\asymp\max\left\{\frac{L}{m}R,\; (\friction^2/\lip) \max(1, \frac{\sqrt{\lip/m}}{\friction}) (m^{-1/2} + (\lip/m)^{3/2}R)\right\}
    \end{equation}
\end{lemma}
\begin{proof}
    Carefully tracking constants in~\citep{schuh2024globalcontractivity}, one can check that $\rcf=0$ if
    \begin{equation}
        |T_t|\gtrsim (\friction^2/\lip) \max(1, \frac{\sqrt{\lip/m}}{\friction}) (m^{-1/2} + (\lip/m)^{3/2}R).
    \end{equation}
    By \ito{}'s lemma for such $T_t$ due to the purely synchronous coupling
    \begin{equation}\label{eq:prof_unif_bound1}
        \begin{aligned}
            \dd T_t 
            &= \begin{bmatrix}
                -\frac{\friction}{2}(X_t-X^\aux{}_t) + \frac{\friction}{2}(\bar P_t-\bar P^\aux{}_t)\\
                - \frac{\friction}{2}(\bar P_t-\bar P^\aux{}_t) + \frac{\friction}{2}(X_t-X^\aux{}_t) - \frac{2}{\friction}(\nabla \pot(X_t)-\nabla \pot(X_t^\aux{}))
            \end{bmatrix} \dd t
        \end{aligned}
    \end{equation}
    and, thus,
    \begin{equation}\label{eq:prof_unif_bound2}
        \begin{aligned}
            \dd |T_t|^2
            =& 2 T_t^\top\begin{bmatrix}
                -\frac{\friction}{2} &  \frac{\friction}{2}\\
                \frac{\friction}{2}& -\frac{\friction}{2}
            \end{bmatrix}
            T_t \dd t
            - \frac{4}{\friction} (\bar P_t-\bar P^\aux{}_t)^\top\Hc_t(X_t-X_t^\aux{}) \dd t \\
        \end{aligned}
    \end{equation}
    where 
    \begin{equation}
        \Hc_t\coloneqq \int_0^1 \nabla^2 \pot (X_t^\aux{} + s(X_t-X_t^\aux{}))\dd s.
    \end{equation}
    Now we make a case distinction. Let $\tilde R>0$ be such that for $|X_t-X_t^\aux{}|\geq\tilde R$ it holds $\Hc_t\succeq \frac{m}{2}I$. This can be ensured for $\tilde R\asymp \frac{L}{m}R$, \cf~\cref{lemma:hessian_bound}. Then, for $|X_t-X_t^\aux{}|\leq\tilde R$ the second row of~\cref{eq:prof_unif_bound1} yields
    \begin{equation}
        \begin{aligned}
            \dd |\bar P_t - \bar P^\aux{}_t|^2 \leq |\bar P_t-\bar P^\aux{}_t|\big(-\friction|\bar P_t-\bar P^\aux{}_t| + \big(\frac{4\lip}{\friction}+ \friction\big)\tilde R\big)
        \end{aligned}
    \end{equation}
    and in particular, $t \mapsto |\bar P_t- \bar P^\aux{}_t|^2$ decreases whenever $|\bar P_t- \bar P^\aux{}_t|\geq \big(\frac{4\lip}{\friction^2} + 1\big)\tilde R$. On the other hand-side, for $|X_t-X_t^\aux{}|\geq\tilde R$,~\cref{eq:prof_unif_bound2} gives
    \begin{equation}
        \begin{aligned}
            \dd |T_t|^2
            \leq& 2 T_t^\top\begin{bmatrix}
                -\frac{\friction}{2} &  \frac{\friction}{2}\\
                \frac{\friction}{2} - \frac{2}{\friction}\Hc_t& -\frac{\friction}{2}
            \end{bmatrix}T_t\dd t\\
            =& 2 T_t^\top\begin{bmatrix}
                -\frac{\friction}{2} &  \frac{\friction}{2} - \frac{1}{\friction}\Hc_t\\
                *& -\frac{\friction}{2}
            \end{bmatrix}T_t\dd t.
        \end{aligned}
    \end{equation}
    where we symmetrized the matrix for the last equality.
    The eigenvalues of the above matrix are the solutions $\lambda$ of the equation
    \begin{equation}
        0 = (\frac{\friction}{2}+\lambda)^2 - (\frac{\friction}{2} - \frac{\mu}{\friction})^2
    \end{equation}
    with $\mu$ an eigenvalue of $\Hc_t$. An explicit computation reveals that if $\friction^2\geq 2 \mu$ the largest eigenvalue is $-\tfrac{\mu}{\friction}$ and since $\mu>\tfrac{m}{2}$ it holds
    \begin{equation}
        \begin{aligned}
            \dd |T_t|^2\leq - \frac{m}{\friction}|T_t|^2\dd t.
        \end{aligned}
    \end{equation}
    In total, since $\friction\geq \sqrt{2\lip}$ it follows that $\ t\mapsto |T_t|^2$ is decreasing when
    \begin{equation}
        |T_t|^2\gtrsim \tilde R^2\gtrsim \frac{L^2}{m^2}R^2.
    \end{equation}
    Combining both cases, noting a.s. path-continuity, accounting for the initialization distance and equivalence of the twisted norm to the Euclidean norm yields the desired result.
\end{proof}

\subsection{Putting everything together}

\begin{theorem}\label{theorem:final_kl_bound}
Let~\cref{ass:pot} be satisfied and $\friction^2 > \frac{2(\lip+m)^2}{m}$. Fix $\delta>0$. There exist dimension free constants $C,\Rc,\chi$ depending only on the problem parameters $\lip,m,R$ and $C$ in addition on $\delta$ such that for $T>\delta$
\begin{equation}
        \begin{aligned}
            \KL(Z_T|Z'_T)
            \leq C (\Rc\vee |z_0-z'_0|) \exp{-\chi T}|z_0-z'_0|.
        \end{aligned}
    \end{equation}
    The constants $\Rc$ is polynomial in the problem parameters. 
\end{theorem}
\begin{proof}
    Note that
    \begin{equation}
        \begin{aligned}
            |T_{\tc}|^2 
            =& |X_{\tc} - X_{\tc}^\aux{}|^2 + |X_{\tc} - X_{\tc}^\aux{} + \frac{2}{\friction_{\tc}}(V_{\tc} - V_{\tc}^\aux{})|^2\\
            \lesssim& |X_{\tc} - X_{\tc}^\aux{}|^2 + |X_{\tc} - X_{\tc}^\aux{} + \frac{2}{\friction}(V_{\tc} - V_{\tc}^\aux{})|^2
        \end{aligned}
    \end{equation}
    since $\friction_{\tc}\geq \friction$. Then, combining~\cref{lemma:interm_bound,lemma:schuh_contraction,lemma:uniform_bound} we find
    \begin{equation}
        \begin{aligned}
            \KL(Z_T|Z'_T)
            \lesssim& \bigg(\frac{|\omega|^3(\exp{-4\omega(T-\tc)}-1)}{\friction(\exp{-\omega(T-\tc)}-1)^4} + \frac{|\omega|\friction(\exp{-2\omega(T-\tc)}-1)}{(\exp{-\omega(T-\tc)}-1)^2}\bigg)(1\vee \friction^{-1})^2(\Rc\vee |z_0-z'_0|) C \exp{- \chi \tc}|z_0-z'_0|\\
        \end{aligned}
    \end{equation}
    Using that $1+t\leq \exp{t}$, $T-\tc = \delta$ and overwriting the constant $C$ we obtain
    \begin{equation}
        \begin{aligned}
            \KL(Z_T|Z'_T)
            \lesssim& \big(\frac{1}{\delta^3} + \frac{1}{\delta}\big)C(\Rc\vee |z_0-z'_0|)|z_0-z'_0|\exp{- \chi T}.
        \end{aligned}
    \end{equation}
\end{proof}

\begin{remark}
    Note that the rate $\chi$ is the exact same convergence rate of~\citep{schuh2024globalcontractivity}. 
    While $\Rc$ is only polynomial, $C,\chi$ deteriorate exponentially in the problem parameters. All constants are independent of the dimension $d$.
\end{remark}

\subsection{Harnack inequalities}\label{sec:harnack}
It is well known that reverse transportation inequalities, such as in \cref{theorem:final_kl_bound}, are dual to Harnack inequalities. We give a proof for completeness.
\begin{lemma}[Logarithmic Harnack Inequality]
    Let $P_t$ be the semigroup of \cref{eq:uld} and $z,z' \in \R^{2d}$. Then \cref{theorem:final_kl_bound} is equivalent to 
    \begin{equation}\label{eq:log_hranack}
        P_t \log f(z) \leq C (\Rc\vee |z-z'|) \exp{-\chi t}|z-z'| + \log P_
        t f(z'),
    \end{equation}
    for all measurable, bounded, and strictly positive $f$.
\end{lemma}
\begin{proof}
    \enquote{$\Rightarrow$}:
    By the Donsker--Varadhan theorem~\cite[Theorem 1.5.7]{Chewi26Book} we have
    \begin{equation}\label{eq:donsker_varadhan}
        \E_\mu[g] \leq \KL(\mu \mid \nu) + \log \E_\nu[\exp{g}],
    \end{equation}
    for probability measures $\mu, \nu$ and any bounded and measurable test function $g$. Now let $z,z' \in \R^{2d}$ and set $\mu = \delta_z P_t$, $\nu = \delta_{z'} P_t$ and set $g = \log f$ for smooth, bounded $f$ satisfying $f\geq \epsilon>0$ for some $\epsilon$.
    Inserting into \cref{eq:donsker_varadhan} and using \cref{theorem:final_kl_bound} to bound the KL term yields
    \begin{align}
        P_t \log f(z) &= \int \log f \dd(\delta_z P_t) 
        \leq \KL(\delta_z P_t\mid \delta_{z'} P_t) + \log \int f \dd (\delta_{z'} P_t)  \\
        &\leq \KL(\delta_z P_t\mid \delta_{z'} P_t) + \log P_t f(z').
    \end{align}
    The result then follows by inserting the upper bound of the $\KL$ from~\cref{theorem:final_kl_bound}.
    For arbitrary bounded, measurable and strictly positive $f$, we may approximate $f$ by $f_n = f \vee \frac{1}{n}$. Then, by monotonicity we have
    \begin{align}
        P_t \log f(z)\leq P_t \log f_n(z)
        &\leq \KL(\delta_z P_t\mid \delta_{z'} P_t) + \log \int f_n \dd (\delta_{z'} P_t)  \\
        &\to \KL(\delta_z P_t\mid \delta_{z'} P_t) + \log P_t f(z').
    \end{align}
    where the convergence follows from the monotone convergence theorem.
    \enquote{$\Leftarrow$}:
    Inserting $f = \exp{g}$ with $g$ bounded and measurable into \cref{eq:log_hranack} yields
    \begin{equation}
        \E_{\delta_z P_t}[g] - \log E_{\delta_{z'}P_t}[\exp{g}] \leq C (\Rc\vee |z-z'|) \exp{-\chi t}|z-z'|.
    \end{equation}
    Taking the supremum over $g$ yields the result again by Donsker--Varadhan.
\end{proof}

\section{Conclusion, limitations and future work}\label{sec:conclusion_futurework}
We provided a proof for a long-term reverse transportation, respectively, Harnack inequality for the underdamped Langevin dynamics in the non-logconcave setting. The proof relies on a coupling construction with an interpolating auxiliary process. In particular, we retain the Wasserstein-1 convergence rates for the underdamped Langevin dynamics from~\citep{schuh2024globalcontractivity}.

\paragraph{Limitations and Future work}
An interesting direction for future work is the adaptation of the results for convergence in R\'enyi divergence. As pointed out in~\citep{lu2026longtimereversetransportationinequalities} in this case it seems necessary to rely on a shift construction with bounded shifts. Moreover, we plan to combine our results with \citep{nagler2026annealedunderdampedlangevindynamics} to obtain convergence results for annealed underdamped Langevin dynamics in the non-convex setting.

\section*{Acknowledgements}
This research was funded in whole or in part by the Austrian Science Fund (FWF) [\href{https://www.fwf.ac.at/forschungsradar/10.55776/PAT4839425}{10.55776/PAT4839425}, \href{https://www.mr-dynamo.at}{10.55776/F100800}]. For open access purposes, the author has applied a CC BY public copyright license to any author-accepted manuscript version arising from this submission.

\bibliography{iclr2027_conference}

@book{karatzas2014brownian,
	title        = {Brownian motion and stochastic calculus},
	author       = {Karatzas, Ioannis and Shreve, Steven},
	year         = 2014,
	publisher    = {springer}
}

@article{leimkuhler2013rational,
	title        = {Rational Construction of Stochastic Numerical Methods for Molecular Sampling},
	author       = {Leimkuhler, Benedict and Matthews, Charles},
	journal      = {Applied Mathematics Research eXpress},
	volume       = 2013,
	number       = 1,
	pages        = {34--56},
	year         = 2012,
	month        = {06}
}

@article{leimkuhler2024contraction,
	title        = {Contraction and Convergence Rates for Discretized Kinetic Langevin Dynamics},
	author       = {Leimkuhler, Benedict J. and Paulin, Daniel and Whalley, Peter A.},
	journal      = {SIAM Journal on Numerical Analysis},
	volume       = 62,
	number       = 3,
	pages        = {1226--1258},
	year         = 2024
}

@inproceedings{cheng2018underdamped,
	title        = {Underdamped Langevin MCMC: A non-asymptotic analysis},
	author       = {Cheng, Xiang and Chatterji, Niladri S and Bartlett, Peter L and Jordan, Michael I},
	booktitle    = {Conference on learning theory},
	pages        = {300--323},
	year         = 2018,
	organization = {PMLR}
}

@book{Chewi26Book,
	title        = {Log-concave sampling},
	author       = {Chewi, Sinho},
	publisher    = {Forthcoming},
	year         = 2026,
	note         = {Available online at \url{https://chewisinho.github.io/}}
}

@article{eberle2026non_rev_lift,
	title        = {Non-reversible lifts of reversible diffusion processes and relaxation times},
	author       = {Eberle, Andreas and L{\"o}rler, Francis},
	journal      = {Probability Theory and Related Fields},
	year         = 2026,
	month        = {Feb},
	day          = {01},
	volume       = 194,
	number       = 1,
	pages        = {173--203},
	issn         = {1432-2064}
}

@article{eberle2017coupling_kinetic,
	title        = {Couplings and quantitative contraction rates for Langevin dynamics},
	author       = {Andreas Eberle and Arnaud Guillin and Raphael Zimmer},
	journal      = {The Annals of Probability},
	year         = 2019
}

@inproceedings{altschuler2026shiftedcompositionivballistic,
  title={Shifted composition IV: toward ballistic acceleration for log-concave sampling},
  author={Altschuler, Jason M and Chewi, Sinho and Zhang, Matthew S},
  booktitle={Proceedings of the 58th Annual ACM Symposium on Theory of Computing},
  pages={1739--1750},
  year={2026}
}

@article{schuh2024globalcontractivity,
	title        = {{Global contractivity for Langevin dynamics with distribution-dependent forces and uniform in time propagation of chaos}},
	author       = {Katharina Schuh},
	volume       = 60,
	journal      = {Annales de l'Institut Henri Poincaré, Probabilités et Statistiques},
	number       = 2,
	publisher    = {Institut Henri Poincaré},
	pages        = {753 -- 789},
	year         = 2024
}

@article{Dalalyan2018OnSF,
	title        = {On sampling from a log-concave density using kinetic Langevin diffusions},
	author       = {Arnak S. Dalalyan and Lionel Riou-Durand},
	journal      = {ArXiv},
	year         = 2018,
	volume       = {abs/1807.09382}
}

@misc{lu2026longtimereversetransportationinequalities,
	title        = {Long-time reverse transportation inequalities for non-globally-dissipative Langevin dynamics},
	author       = {Jianfeng Lu and Yuliang Wang},
	year         = 2026,
	eprint       = {2512.18598},
	archiveprefix = {arXiv},
	primaryclass = {math.PR}
}

@article{altschuler2025shifted_comp_1,
	title        = {Shifted Composition I: Harnack and Reverse Transport Inequalities},
	author       = {Altschuler, Jason M. and Chewi, Sinho},
	year         = 2025,
	issue_date   = {Jan. 2025},
	publisher    = {IEEE Press},
	volume       = 71,
	number       = 1,
	issn         = {0018-9448},
	journal      = {IEEE Trans. Inf. Theor.},
	month        = jan,
	pages        = {90–113},
	numpages     = 24
}

@article{arnaudon2006harnackinequalityheatkernel,
	title        = {Harnack inequality and heat kernel estimates on manifolds with curvature unbounded below},
	journal      = {Bulletin des Sciences Mathématiques},
	volume       = 130,
	number       = 3,
	pages        = {223--233},
	year         = 2006,
	issn         = {0007-4497},
	author       = {Marc Arnaudon and Anton Thalmaier and Feng-Yu Wang}
}

@misc{wang2010couplingapplications,
	title        = {Coupling and Applications},
	author       = {Feng-Yu Wang},
	year         = 2010,
	eprint       = {1012.5687},
	archiveprefix = {arXiv},
	primaryclass = {math.PR}
}

@article{shen2019randomized,
  title={The randomized midpoint method for log-concave sampling},
  author={Shen, Ruoqi and Lee, Yin Tat},
  journal={Advances in Neural Information Processing Systems},
  volume={32},
  year={2019}
}

@article{villani2006hypocoercivity,
  title={Hypocoercivity},
  author={Villani, C{\'e}dric},
  journal={arXiv preprint math/0609050},
  year={2006}
}

@misc{nagler2026annealedunderdampedlangevindynamics,
      title={Annealed Underdamped Langevin Dynamics}, 
      author={Laurenz Nagler and Alexander Falk and Andreas Habring},
      year={2026},
      eprint={2609.32446},
      archivePrefix={arXiv},
      primaryClass={math.PR},
      url={https://arxiv.org/abs/2609.32446}, 
}

@misc{guillin2012degeneratefokkerplanckequations,
      title={Degenerate Fokker-Planck Equations : Bismut Formula, Gradient Estimate and Harnack Inequality}, 
      author={Arnaud Guillin and Feng-Yu Wang},
      year={2012},
      eprint={1103.2817},
      archivePrefix={arXiv},
      primaryClass={math.PR},
      url={https://arxiv.org/abs/1103.2817}, 
}

@article{wang1997logarithmic,
  title={Logarithmic Sobolev inequalities on noncompact Riemannian manifolds},
  author={Wang, Feng-Yu},
  journal={Probability theory and related fields},
  volume={109},
  number={3},
  pages={417--424},
  year={1997},
  publisher={Springer}
}

@book{wang2013harnack,
  title={Harnack inequalities for stochastic partial differential equations},
  author={Wang, Feng-Yu},
  volume={1332},
  year={2013},
  publisher={Springer}
}

@article{wang2010harnack,
  title={Harnack inequalities on manifolds with boundary and applications},
  author={Wang, Feng-Yu},
  journal={Journal de math{\'e}matiques pures et appliqu{\'e}es},
  volume={94},
  number={3},
  pages={304--321},
  year={2010},
  publisher={Elsevier}
}

@article{wang2011harnack,
  title={Harnack inequalities for functional SDEs with multiplicative noise and applications},
  author={Wang, Feng-Yu and Yuan, Chenggui},
  journal={Stochastic processes and their applications},
  volume={121},
  number={11},
  pages={2692--2710},
  year={2011},
  publisher={Elsevier}
}
\bibliographystyle{abbrvnat}

\appendix

\section{Rigorous arguments for the application of Girsanov's theorem}\label{sec:rigorous_girsanov}
Let \((\Omega, \Ac, (\Fc_t)_{t \geq 0}, \P)\) denote a filtered probability space satisfying the \emph{usual assumptions}. For \(z_0, z_0' \in \R^{2d}\) and $\tc\in [0,T]$ we define the coupled \gls{sde},
\begin{equation}\label{eq:coupling_construction}
    \begin{aligned}
    &\begin{cases}
        \dd X_t &= V_t \dd t\\
        \dd V_t &= \{-\friction V_t - \nabla \pot (X_t) \}\,\dd t + \sqrt{2 \friction}\dd W_t\\
        \dd X_t^{\aux{}} &= V_t^{\aux{}}\dd t,\\
        \dd V_t^{\aux{}} &=
        \begin{aligned}[t]
             &\{-\friction V_t^{\aux{}}
            - \nabla \pot (X^\aux{}_t) \\
            & + \eta_t^x(X_t - X_t^{\aux{}})
            + \eta_t^v(V_t - V_t^{\aux{}})\}\dd t
            + \sqrt{2 \friction}\,\dd \bar W_t.
        \end{aligned}
    \end{cases}\\
    &Z_0 = z_0,\; 
    Z_0^\aux{}= z_0'
    \end{aligned}\tag{AUX}
\end{equation}
where the shift parameters $\eta^x_t, \eta^v_t$ are to be defined later. For the coupling of the Brownian motions $(W_t)_t$, $(\bar W_t)_t$ we refer to~\cref{sec:coupling}, respectively,~\citep{schuh2024globalcontractivity}.

Next, define the measure $\Q$ via $N_t\coloneqq -\eta_t^x(X_t - X_t^{\aux{}}) - \eta_t^v(V_t - V_t^{\aux{}})$ and 
\begin{equation}
    \frac{\dd \Q}{\dd \P} = \exp{\frac{1}{\sqrt{2\friction}}\int_0^t N_s \dd \bar W_s - \frac{1}{4\friction}\int_0^t |N_s|^2\dd s}.
\end{equation}
respectively,
\begin{equation}
    \frac{\dd \P}{\dd \Q} = \exp{\frac{1}{\sqrt{2\friction}}\int_0^t - N_s \dd \bar W_s + \frac{1}{4\friction}\int_0^t |N_s|^2\dd s}\eqqcolon M_t((Z_s,Z^\aux{}_s)_{s\leq t})
\end{equation}
where $M_t$ is a function purely of $(Z_s,Z^\aux{}_s)_{s\leq t}$ since also the Brownian motion $(\bar W_t)_t$ can be expressed in terms of this process.
By Girsanov's theorem $(B_t)_t$ defined via $B_0=0$ and 
\begin{equation}
    \sqrt{2 \friction}\dd B_t = \{\eta_t^x(X_t - X_t^{\aux{}})
    + \eta_t^v(V_t - V_t^{\aux{}})\}\dd t
    + \sqrt{2 \friction}\,\dd \bar W_t.
\end{equation}
is a $\Q$-Brownian motion and~\cref{eq:coupling_construction} satisfies $\Q$-a.s.
\begin{equation}\label{eq:coupling_construction}
    \begin{aligned}
    &\begin{cases}
        &\dd X_t^{\aux{}} = V_t^{\aux{}}\dd t \\
        &\dd V_t^{\aux{}} =
        \begin{aligned}[t]
             &\{-\friction V_t^{\aux{}}
            - \nabla \pot(X_t^{\aux{}})\}\dd t
            + \sqrt{2 \friction}\,\dd B_t.
        \end{aligned}\\
    &(X_0^\aux{},V_0^\aux{}) = (x_0',v_0')
    \end{cases}
    \end{aligned}\tag{AUX}
\end{equation}
The $\Q$-distribution of the first components $(X_t,V_t)_t$ in~\cref{eq:coupling_construction} is not directly tractable as it is unclear how the Brownian motion $(W_t)_t$ behaves under the change of measure, but it is also irrelevant for our purposes.
If we denote as $\mu^{z_0}_{\sh{}}$ the path measure of~\cref{eq:coupling_construction} in $\P$ and as $\mu^{z_0}_{\ud}$ the one in $\Q$ we have
\begin{equation}
    \frac{\dd \mu^{z_0}_{\sh{}}}{\mu^{z_0}_{\ud}}( (z_s,z^\aux{}_s)_{s\leq t}) = M_t((z_s,z^\aux{}_s)_{s\leq t})
\end{equation}
It follows
\begin{equation}
    \begin{aligned}
        \KL(Z^\aux{}|Z)
        =& \E_{\P}[\log \frac{\dd \mu^{z_0}_{\sh{}}}{\mu^{z_0}_{\ud}}((Z_s,Z^\aux{}_s)_{s\leq t})]\\
        =& \int_0^t \E_{\P}[\frac{1}{4\friction}|\eta_t^x(X_t - X_t^{\aux{}}) + \eta_t^v(V_t - V_t^{\aux{}})|^2]\dd s.
    \end{aligned}
\end{equation}

\begin{remark}
    The application of Girsanov's theorem can be justified using standard localization techniques.
\end{remark}

\section{Helper Results}

\begin{lemma}\label{lemma:hessian_bound}
Let~\cref{ass:pot} hold. Define
\begin{equation}
    \Hc(x,y) \coloneqq \int_0^1 \nabla^2 U(y+t(x-y)) \dd t.
\end{equation}
Then $|x-y|\geq R'$ implies $\Hc(x,y)\succeq m'\coloneqq\frac{R'-2R}{R'}m - L \frac{2R}{R'}$. In particular, if $|x-y| \geq 4(1+\tfrac{\lip}{m})R$, then $\Hc(x,y)\succeq \tfrac{m}{2}$.
\end{lemma}
\begin{proof}
    Let $w\in\R^d\setminus\{0\}$. Indeed, by~\cref{ass:pot}
    \begin{equation}
        \begin{aligned}
            w^\top \int_0^1 \nabla^2 U(y+t(x-y))\dd t w
            \geq& |w|^2\int_0^1 \1_{|y+t(x-y)|\geq R}m - \1_{|y+t(x-y)|<R}L \dd t\\
            \geq& ((1-\lambda)m - \lambda L)|w|^2\\
        \end{aligned}
    \end{equation}
    where we define $\lambda\in [0,1]$ is defined as 
    \begin{equation}
        \lambda = |\{t\in[0,1]\;|\;|y+t(x-y)|<R\}|
    \end{equation}
    with $|\cdot|$ the Lebesgue measure.
    One can easily verify that for $|x-y|\geq R'$, $\lambda \leq 1 \wedge \frac{2R}{R'}$ implying the desired result.
\end{proof}

\end{document}